\documentclass{article}[12pt]
\usepackage{amsmath}
\usepackage{amsthm}
\usepackage{amsfonts}
\usepackage{amssymb}
\usepackage{amstext}
\usepackage{amsopn}
\usepackage{latexsym}
\newtheorem{thm}{Theorem}

\newtheorem{cor}{Corollary}
\newtheorem{lemma}{Lemma}

\begin{document}
\title{Doubling measures on uniform Cantor sets\thanks{Supported
by the NSFC (Nos. 11271114 and 10971056).}}

\author{Chun WEI,
\quad  Shengyou WEN,\thanks{Corresponding author. E-mail:
sywen\_65@163.com}\quad and Zhixiong WEN}

\maketitle

\begin{abstract}
We obtain a complete description for a probability measure to be doubling on an arbitrarily
given uniform Cantor set. The question of which doubling measures on such a Cantor set can be
extended to a doubling measure on $[0,1]$ is also considered.

\medskip


\end{abstract}

\section{Introduction}

It is well known that every complete doubling metric space carries a doubling measure
\cite{LS,KRS,VK} and that doubling measures are rich
when the underlying doubling metric space has no dense isolated points \cite{KW}.
In the Euclidean $n$-space, a measure, defined by integrating an $A_\infty$
weight, is doubling \cite{G}. These measures form a subclass of absolutely
continuous doubling measures on $\mathbb{R}^n$.
Examples of singular doubling measures can be obtained from self-similar measures
on the unit square \cite{Fa,Y}.
The existence of singular doubling measures have also been studied extensively
in general metric spaces \cite{KW,Tu,Wu}. Recall that a Borel regular measure $\mu$ on a
metric space $X$ is doubling, if there is a constant $C\geq 1$ such that
\begin{equation}\label{eq1}
0<\mu(B(x,2r))\leq C\mu(B(x,r))<\infty\,\,\mbox{for all $x\in X$ and $r>0$.}
\end{equation}

Doubling measures naturally arise in the study of the gradients of
convex functions on the Euclidean $n$-space. See \cite{KM,KMW}, where it is showed that,
in many cases, the gradient of
a convex function is a quasisymmetric self-homeomorphism of $\mathbb{R}^n$, by which
the pullback measure of the Lebesgue measure is doubling. Doubling measures can be
applied to classify
subsets of a metric space, in which, for example,
a subset is called very thin, if it is null for all doubling measures on the metric space.
This leads to the study on fat and thin sets \cite{B,hww,ORS}.
Doubling measures can also be used to construct quasisymmetric deformations of the
underlying metric space, so that the objective space is in some sense better than the original space,
for example, it can be proved that, if a metric space is uniformly perfect and carries a doubling
measure, then it is quasisymmetrically
equivalent to an Ahlfors regular space \cite{H1}.

Though examples of doubling measures are rich, a complete description for doubling measures
on a given metric space is usually difficult. It is clear that a measure on the real line
$\mathbb{R}$ is doubling if and only if it is a pullback measure of a quasisymmetric self-homeomorphism
of $\mathbb{R}$. However, this is not true for higher dimensional Euclidean
space. In the present paper, we study doubling measures on uniform Cantor sets.
We shall give a sufficient and necessary condition
for the doubling property of a measure on an arbitrarily given uniform Cantor set.
The question of which doubling measure on such a Cantor set can be extended to a doubling
measure on $[0,1]$ is also discussed. For the related papers on Cantor sets and doubling measures
on $[0,1]$, we refer to \cite{CS,PW,SW,WW}.

We begin with uniform Cantor sets. Let
$\mathbf{n}=\{n_k\}_{k=1}^\infty$ be a sequence of
integers, where $n_k\geq2$. Let $\mathbf{c}=\{c_k\}_{k=1}^\infty$ be a sequence of real
numbers in $(0,1)$ such that $(n_k-1)c_k<1$ for all $k\geq1$. The
uniform Cantor set $E(\mathbf{n},\mathbf{c})$ is defined by
\begin{equation}\label{eq2}
E(\mathbf{n},\mathbf{c})=\bigcap_{k=0}^\infty E_k,
\end{equation}
where $\{E_k\}$ is a sequence of nested compact sets in $[0,1]$,
$E_0=[0,1]$, and $E_k$ is obtained by deleting from every
component $I$ of $E_{k-1}$ $(n_k-1)$ open intervals of equal length
$c_k|I|$, such that the remaining $n_k$ closed intervals are of
equal length. A component of $E_k$ will be called a component of
level $k$ and a component of $E_{k-1}\setminus E_k$ will be called
a gap of level $k$. Denote by $N_k$, $\delta_k$, and
$\varepsilon_k$, respectively, the number of components of level
$k$, the length of a component of level $k$, and the length of a
gap of level $k$. Then, from the definition
\begin{equation}\label{eq3}
N_k=\prod_{i=1}^kn_i,\ \
\delta_k=\prod_{i=1}^k\frac{1-(n_i-1)c_i}{n_i},\ \
\varepsilon_k=c_k\delta_{k-1}.
\end{equation}
Note that the uniform Cantor set $E(\mathbf{n},\mathbf{c})$
is of Lebesgue measure
$$\mathcal{L}(E(\mathbf{n},\mathbf{c}))=\prod_{i=1}^\infty(1-(n_i-1)c_i).$$
Therefore
\begin{equation}\label{eq4}
\mathcal{L}(E(\mathbf{n},\mathbf{c}))>0\Leftrightarrow\{n_kc_k\}\in\ell^1,
\end{equation}
where $\ell^p$ denotes the set of real sequences $\{a_k\}_{k=1}^\infty$ with $\sum_{k=1}^\infty
|a_k|^p<\infty$.

For a uniform Cantor set
$E(\mathbf{n},\mathbf{c})$ we denote by $\mathcal{I}_k$ the family of its components of level $k$
and by $\mathcal{G}_k$ the family of its gaps of level $k$. Then $\sharp\mathcal{I}_k=N_k$
and $\sharp\mathcal{G}_k=N_{k-1}(n_k-1)$, where $\sharp$ denotes the cardinality.
Let $\mathcal{I}=\cup_{k=1}^\infty\mathcal{I}_k$ and $\mathcal{G}=\cup_{k=1}^\infty\mathcal{G}_k$.
To label these components and gaps, we relate a word to each of them as
follows: Let
\begin{equation}\label{eq5}
W_k=\{i_1i_2\cdots i_k : 1\leq i_j\leq n_j, 1\leq j\leq
k\}\,\,\mbox{and}\,\, W=\cup_{k=1}^\infty W_k.
\end{equation}
A member in $W_k$ is called
a word of length $k$ and a member in $W$ is called a word of finite
length. Let the words in $W$ be ordered in the lexicographic order.
Then the first $n_1$ words are $1,2,\cdots,n_1$. We denote, from left to right, by $I_1,I_2,\cdots,I_{n_1}$
the $n_1$ components in $\mathcal{I}_1$ and by
$G_1,G_2,\cdots,G_{n_1-1}$ the $n_1-1$ gaps in $\mathcal{G}_1$. Inductively,
supposing that $I_w$ has been defined for a word $w\in W_{k-1}$, we denote, from left to right, by
$I_{w1},I_{w2},\cdots,I_{wn_k}$ the $n_k$ members of
$\mathcal{I}_k$ in $I_w$ and by $G_{w1},G_{w2},\cdots,G_{w(n_k-1)}$ the
$n_k-1$ members of $\mathcal{G}_k$ in $I_w$. With the above notation, once a word $w\in W_k$ is
specified, the relative position of the component $I_w$ of level $k$ is
determined. For example, the minimal word $11\cdots 1$ of length $k$ corresponds to the
leftmost component in $\mathcal{I}_k$ and the maximal word $n_1n_2\cdots n_k$ of length $k$
corresponds to the rightmost component in $\mathcal{I}_k$.
We say that two components $I_{w}$ and $I_{u}$ are adjacent, if $w,u\in W$ are adjacent in the sense
of the lexicographic order. For example, given $w\in W_{k-1}$ and $1\leq i<n_k$,
the maximal word $win_{k+1}\cdots n_{k+t}$ of length $k+t$ with prefix $wi$ and the minimal word $w(i+1)1\cdots 1$ of
length $k+t$ with prefix $w(i+1)$ are adjacent. They correspond to two adjacent components in $\mathcal{I}_{k+t}$. Note that
$$E(\mathbf{n},\mathbf{c})=\bigcap_{k=1}^\infty
\bigcup_{\omega\in W_k}I_\omega$$ and that
$$I_\omega =I_{w1}\cup G_{w1}\cup I_{w2}\cup G_{w2}\cup\cdots\cup I_{w(n_k-1)}\cup G_{w(n_k-1)}\cup I_{wn_k}$$
for all $w\in W_{k-1}$, where, if $k=1$, then $I_w=[0,1]$ and $wi=i$.

We also need some terminology on vectors. A vector $P=(p_1,p_2,\cdots,p_k)$ is called a positive probability,
if its components are all positive and $\sum_{i=1}^kp_i=1$.
A positive probability vector $P=(p_1,p_2,\cdots,p_k)$ is called $C$-uniform,
where $C\geq 1$ is a constant, if
\begin{equation}\label{jia}
C^{-1}\leq \frac{
p_{i+1}+\cdots+p_{i+l}}{p_{j+1}+\cdots+p_{j+l}}
\leq C\,\,\mbox{whenever $0\leq i\leq j\leq i+l\leq j+l\leq k$.}
\end{equation}
The condition means that two of
sums of consecutive $l$ components of $P$ are comparable, if they are adjacent.
Clearly, the vector $P$ is $1$-uniform if and only if $p_1=p_2=\cdots=p_k$.
We say that $P$ is $(C,s)$-uniform, where $1\leq s\leq k$, if
the condition (\ref{jia}) holds only for $s\leq l\leq k$. Thus, the
$(C,1)$-uniformness and the $C$-uniformness are the same. We say that a
sequence $\mathbb{P}=\{P_k\}$ of positive probability vectors
is $C$-uniform, if each $P_k$ is $C$-uniform. We say that $\mathbb{P}$ is ultimately $C$-uniform,
if there is an integer $k_0\geq1$ such that $P_k$ is $C$-uniform for each $k\geq k_0$.

Let $E(\mathbf{n},\mathbf{c})$ be a uniform Cantor set.
Let $\mathbb{P}=\{P_k\}_{k=1}^\infty$ be a sequence of positive probability vectors, where
$$P_k=(p_{k,1},p_{k,2},\cdots,p_{k, n_k})$$
has $n_k$ components.
Such a sequence of vectors will be called $\mathbf{n}$-matching. Given an $\mathbf{n}$-matching
sequence $\mathbb{P}$ of vectors, there is a unique probability measure on
$E(\mathbf{n},\mathbf{c})$, denoted by $\mu_\mathbb{P}$, satisfying
\begin{equation}\label{eq6}
\mu_\mathbb{P}(I_{wi})=p_{k,i}\mu_\mathbb{P}(I_w)\mbox{ for all $w\in W_{k-1}$ and $1\leq i\leq n_k$.}
\end{equation}
The positivity assumption of probability vectors is reasonable. Otherwise, the corresponding measure would support on
a proper compact subset of $E(\mathbf{n},\mathbf{c})$.

The following data of the uniform Cantor set $E(\mathbf{n},\mathbf{c})$ are crucial in stating our results.
Denote
\begin{equation}\label{eq7}
\Lambda=\{k\geq1: \varepsilon_k<\delta_k\}.
\end{equation}
Thus, a gap in $\mathcal{G}_k$ is shorter than a component in $\mathcal{I}_k$ at every level $k\in\Lambda$.
For each $k\in\Lambda$ let $m_k$ be the unique positive integer such
that
\begin{equation}\label{eq8}
\delta_{k+m_k}\leq\varepsilon_k<\delta_{k+m_k-1}.
\end{equation}
and let $s_k$ be the unique positive
integer such that
\begin{equation}\label{eq9}
s_k\delta_{k+m_k}+(s_k-1)\varepsilon_{k+m_k}\leq\varepsilon_k<(s_k+1)\delta_{k+m_k}+s_k\varepsilon_{k+m_k}.
\end{equation}
Write $E(\mathbf{n},\mathbf{c},\Lambda, \{m_k\}_{k\in\Lambda}, \{s_k\}_{k\in\Lambda})$
for $E(\mathbf{n},\mathbf{c})$, when these related data are emphasized.

\begin{thm}\label{1}
Let $E=E(\mathbf{n},\mathbf{c},\Lambda, \{m_k\}_{k\in\Lambda}, \{s_k\}_{k\in\Lambda})$ be a uniform
Cantor set. Let $\mu_\mathbb{P}$ be a measure determined by an $\mathbf{n}$-matching sequence
$\mathbb{P}$ of positive probability vectors. Then $\mu_\mathbb{P}$ is doubling on $E$ if and only
if there is a constant $C>1$ such that $\mathbb{P}$ is $C$-uniform and satisfies, for each $k\in\Lambda$,
the following conditions:

$(a)$ $(P_{k+t},P_{k+t})$ is $C$-uniform and
\begin{equation}\label{eq11}
C^{-1}\prod_{j=1}^tp_{k+j,1}\leq\prod_{j=1}^t{p_{k+j, n_{k+j}}}\leq C\prod_{j=1}^tp_{k+j,1}\,\,\mbox{for all\,\, $1\leq t<m_k$.}
\end{equation}

$(b)$ $(P_{k+m_k},P_{k+m_k})$ is $(C,s_k)$-uniform.

\noindent Hereafter, if $P=(p_1,p_2,\cdots,p_j)$, then $(P,P)=(p_{1},p_2,\cdots,p_{j},p_{1},p_2,\cdots,p_{j})$.
\end{thm}

From theorem \ref{1} some easier conditions for
the doubling property of the measure $\mu_{\mathbb{P}}$ on $E$ can be formulated.

\begin{cor}\label{c1}
If $\sup _{k\geq1}n_k<\infty$, then the measure $\mu_\mathbb{P}$ is
doubling on $E$ if and only
if there is a constant $C>1$ such that

(a)\, $C^{-1}p_{k,i}\leq p_{k,i+1}\leq Cp_{k,i}$ for all $1\leq i<n_k$ and for all $k$, and

(b)\, the condition (\ref{eq11}) is satisfied for all $1\leq t<m_k$ and for all $k\in\Lambda$.
\end{cor}

\begin{cor}\label{c2}
If $\sup_{k\in\Lambda}m_k<\infty$, then the measure $\mu_\mathbb{P}$
is doubling on $E$ if and only
if there is a constant $C>1$ such that

(a)\, $\mathbb{P}$ is $C$-uniform, and

(b)\, for every integer $k\in\Lambda$, $(P_{k+t},P_{k+t})$ is $C$-uniform for all $1\leq t<m_k$
and $(P_{k+m_k},P_{k+m_k})$ is $(C,s_k)$-uniform.
\end{cor}

\begin{cor}\label{c0}
If $\mathbb{P}$ is ultimately $1$-uniform,
then $\mu_{\mathbb{P}}$ is doubling on $E$.
\end{cor}

Theorem \ref{1} can be generalized. Let $\mu$ be a Borel probability measure on
$E$ with $\mu(I_w)>0$ for all $w\in
W$. For each $k\geq 1$ and for each $w\in W_{k-1}$ let
\begin{equation}\label{eq66}
p_{wi}=\frac{\mu(I_{wi})}{\mu(I_w)},\,\,
1\leq i\leq n_k.
\end{equation}
Then there is a one-to-one
correspondence between the set of Borel probability measures $\mu$
on $E$, with $\mu(I_w)>0$ for all $w\in W$,
and the set of sequences $\{p_w\}_{w\in W}$, with
$\sum_{i=1}^{n_k}p_{wi}=1$ and $0<p_{wi}<1$ for all $w\in W_{k-1}$ and $k\geq 1$. Write
\begin{equation}\label{cr}
P_w=(p_{w1},p_{w2},\cdots,p_{wn_k})
\end{equation}
for each $w\in W_{k-1}$. Then $\mu$ is a measure considered in theorem \ref{1}, provided that
for every fixed integer $k\geq 1$, the vectors $P_w$ are the same for all $w\in W_{k-1}$.
Write $wiu_t$ for the maximal word in $W_{k+t}$ with prefix $wi$
and $wi1^t$ the minimal word in $W_{k+t}$ with prefix $wi$. Our next theorem gives a complete
description of doubling measures on a uniform Cantor set.

\begin{thm}\label{t11}
Let $E=E(\mathbf{n},\mathbf{c},\Lambda, \{m_k\}_{k\in\Lambda}, \{s_k\}_{k\in\Lambda})$ be a
uniform Cantor set. Let $\mu$ be a Borel probability measure on $E$ with $\mu(I_w)>0$ for all
$w\in W$. Let $\{p_w\}_{w\in W}$ be the corresponding sequence.
Then $\mu$ is doubling on $E$ if and only
if there is a constant $C>1$ such that $P_w$ is $C$-uniform for all ${w\in W}$ and satisfies,
for each $k\in\Lambda$, the following conditions:

$(a)$ $(P_{wiu_{t-1}},P_{w(i+1)1^{t-1}})$ is $C$-uniform and
\begin{equation}\label{eq13}
C^{-1}\prod_{j=1}^{t}p_{wiu_j}\leq\prod_{j=1}^{t}{p_{w(i+1)1^j}}\leq C\prod_{j=1}^{t}p_{wiu_j}
\end{equation}
for all $w\in W_{k-1}$, $1\leq i<n_k$, and $1\leq t<m_k$.

$(b)$ $(P_{wiu_{m_k-1}},P_{w(i+1)1^{m_k-1}})$ is $(C,s_k)$-uniform.
\end{thm}

For a uniform Cantor set $E(\mathbf{n},\mathbf{c})$ it is known that
$\nu(E(\mathbf{n},\mathbf{c}))=0$ for all doubling measures $\nu$ on $[0,1]$
if and only if
$\{n_kc_k\}_{k=1}^\infty\not\in\cup_{q\geq1}\ell^q$ \cite{hww}. As a consequence,
if $E(\mathbf{n},\mathbf{c})$ carries a doubling measure which can be extended to a
doubling measure on $[0,1]$, then $\{n_kc_k\}_{k=1}^\infty\in\cup_{q\geq1}\ell^q$.
Our main result on the extension question of doubling measures on a uniform Cantor
set is the following theorem.

\begin{thm}\label{t2}
Let $E=E(\mathbf{n},\mathbf{c})$ be a uniform Cantor set.
Let $\mu_{\mathbb{P}}$ be a measure determined by an ultimately $1$-uniform $\mathbf{n}$-matching
sequence $\mathbb{P}$ of positive probability vectors. Then $\mu_{\mathbb{P}}$ can be
extended to a doubling measure on $[0,1]$ if and only if $\{n_kc_k\}_{k=1}^\infty\in\ell^1$.
\end{thm}

\section{Proof of Theorem \ref{1}}

In what follows
$C$ will denote a constant depending only on the constants in question and it
may be different in every appearance. Write $A\sim B$
for $C^{-1}B\leq A\leq CB$, and $A\preceq B$
for $A\leq CB$, when $A, B$ are quantities. For an interval $I$, denote by
$x^-(I)$ and $x^+(I)$,
respectively, the left and the right endpoints of $I$.

\medskip

\noindent{\bf Proof of the `only if' part.\,}
Let $E=E(\mathbf{n},\mathbf{c},\Lambda, \{m_k\}_{k\in\Lambda}, \{s_k\}_{k\in\Lambda})$ be a uniform
Cantor set. Let $\mu_\mathbb{P}$ be a measure determined by an $\mathbf{n}$-matching sequence
$\mathbb{P}$ of positive probability vectors. Suppose that $\mu_\mathbb{P}$
is doubling on $E$. We are going to show that the sequence $\mathbb{P}$ satisfies
the conditions of theorem \ref{1}. The doubling property of the measure $\mu_\mathbb{P}$ will be used frequently in the following equivalent form:
$$\mu(B(x, r))\sim\mu(B(y,t))\mbox{ for all $x,y\in E$ and $r, s>0$ with $|x-y|\preceq r\sim s$}.$$
Note that, for all $k\geq1$, the endpoints of components of level $k$ of $E$ belong to $E$. The proof consists of the following claims.

\medskip

{\bf{Claim 1.}} $p_{k,2}\leq Cp_{k,1}$ and $p_{k,n_k-1}\leq Cp_{k,n_k}$ for all $k\geq 1$.

\medskip

Proof of Claim 1. Let $w\in W_{k-1}$ be the minimal word of length $k-1$.
Then the corresponding component $I_w$
is the leftmost component of level $k-1$. We easily see that
$$I_{w1}=E_k\cap B(x^-(I_{w1}),\delta_k+\varepsilon_k)\,\,\mbox{and}
\,\, I_{w2}\subseteq B(x^-(I_{w2}),\delta_k+\varepsilon_k).$$
Observing that $x^-(I_{w1}),x^-(I_{w2})\in E$
and $|x^-(I_{w1})-x^-(I_{w2})|=\delta_k+\varepsilon_k$, it follows from the
doubling property of $\mu_\mathbb{P}$ that
$$\mu_\mathbb{P}(I_{w2})\leq\mu_\mathbb{P}(B(x^-(I_{w2}),\delta_k+\varepsilon_k))\leq C\mu_\mathbb{P}(B(x^-(I_{w1}),\delta_k+\varepsilon_k))=C\mu_\mathbb{P}(I_{w1}),$$
giving $p_{k,2}\leq Cp_{k,1}$ by (\ref{eq6}). Similarly, we have $p_{k,n_k-1}\leq Cp_{k,n_k}$
by considering the maximal word $w\in W_{k-1}$ and the rightmost two components of level $k$.
This proves claim 1.

\medskip

{\bf{Claim 2.}} $p_{k,i}\sim p_{k,i+1}$ for all $1\leq i<n_k$ and $k\geq 1$.

\medskip

Proof of Claim 2. Let $k$ be given. When $n_k=2$, claim 2 follows from claim 1 directly.
Next we assume that $n_k>2$. Consider the components $I_{wi}$ and $I_{w(i+1)}$ of level $k$,
where $w\in W_{k-1}$ and $1\leq i< n_k$. Then $G_{wi}$ is the gap of level $k$ between them.
Three possible cases may happen.

\medskip

Case 1. $\varepsilon_k\geq\delta_k$.

\medskip

Since $\varepsilon_k\geq\delta_k$, we easily see that
\begin{equation}\label{peq2}
I_{wi}=E_k\cap B(x^+(I_{wi}),\varepsilon_k)\,\,\mbox{and}\,\, I_{w(i+1)}=E_k\cap B(x^-(I_{w(i+1)}),\varepsilon_k)
\end{equation}
for $1<i<n_k-1$. Observing that $x^+(I_{wi}),x^-(I_{w(i+1)})\in E$ and $|x^+(I_{wi})-x^-(I_{w(i+1)})|=\varepsilon_k$,
we get from doubling property of $\mu_\mathbb{P}$ that
$$\mu_\mathbb{P}(I_{wi})\sim\mu_\mathbb{P}(I_{w(i+1)}),$$
which yields $p_{k,i}\sim p_{k,i+1}$ by (\ref{eq6}). For $i=1$, letting $w$ be the minimal word in $W_{k-1}$,
the relationship (\ref{peq2}) remains to be true, and so
$p_{k,1}\sim p_{k,2}$. For $i=n_k-1$, considering the maximal word in $W_{k-1}$,
we get $p_{k,n_k-1}\sim p_{k,n_k}$ in the same way.

\medskip

Case 2. $\varepsilon_k<\delta_k$ and $n_{k+1}>2$.

\medskip

Consider balls $B(x,r)$ and $B(y,r)$, where
$$x=x^-(G_{wi})-r,\,\, y=x^+(G_{wi})+r,\,\, r=[\frac{n_{k+1}-1}{2}](\delta_{k+1}+\varepsilon_{k+1}).$$
We see that $x$ is an endpoint of a component of level $k+1$ in $I_{wi}$ and that $y$ is
an endpoint of a component of level $k+1$ in $I_{w(i+1)}$, so $x,y\in E$.
Note also that
$$B(x,r)\subseteq I_{wi}\,\,\mbox{and}\,\, B(y,r)\subseteq I_{w(i+1)}.$$
Since $n_{k+1}\geq 3$, we have
$4[\frac{n_{k+1}-1}{2}]\geq n_{k+1}$, so $4r>\delta_k$, giving
$$I_{wi}\subseteq B(x,4r)\,\,\mbox{and}\,\, I_{w(i+1)}\subseteq B(y,4r).$$ Therefore
$$\mu_{\mathbb{P}}(I_{wi})\sim\mu_{\mathbb{P}}(B(x,r))\,\,\mbox{and}\,\,\mu_{\mathbb{P}}(I_{w(i+1)})\sim\mu_{\mathbb{P}}(B(y,r)).$$
Additionally, since $\varepsilon_k<\delta_k$ has been assumed, we have $|x-y|=2r+\varepsilon_k\leq6r$, so
$\mu_{\mathbb{P}}(B(x,r))\sim\mu_{\mathbb{P}}(B(y,r))$ by the doubling property of $\mu_{\mathbb{P}}$.
It then follows that $\mu_{\mathbb{P}}(I_{wi})\sim\mu_{\mathbb{P}}(I_{w(i+1)})$, giving
$p_{k,i}\sim p_{k,i+1}$ by (\ref{eq6}).

\medskip
Case 3. $\varepsilon_k<\delta_k$ and $n_{k+1}=2$.
\medskip

In the case of $\varepsilon_{k+1}\geq\delta_{k+1}$, consider $B(x^-(I_{wi2}),r)$ and $B(x^+(I_{w(i+1)1}),r)$,
where $r=\delta_{k+1}+\frac{\varepsilon_k}{3}$. Since $r<\delta_{k+1}+\frac{\delta_k}{3}\leq\delta_{k+1}+\varepsilon_{k+1}$,
we see that
$$I_{wi2}\subseteq E_{k}\cap B(x^-(I_{wi2}),r)\subseteq I_{wi}$$ and
$$I_{w(i+1)1}\subseteq
E_{k}\cap B(x^+(I_{w(i+1)1}),r)\subseteq I_{w(i+1)}.$$
Observing that
$|x^+(I_{w(i+1)1}-x^-(I_{wi2})|=2\delta_{k+1}+\varepsilon_k\leq3r$,
we get from the doubling property of $\mu_{\mathbb{P}}$ that
$$\mu_{\mathbb{P}}(B(x^-(I_{wi2}),r))\sim\mu_{\mathbb{P}}(B(x^+(I_{w(i+1)1}),r)).$$ Since $n_{k+1}=2$, we have $\mu_{\mathbb{P}}(I_{wi2})\sim\mu_{\mathbb{P}}(I_{wi})$
and $\mu_{\mathbb{P}}(I_{w(i+1)1})\sim\mu_{\mathbb{P}}(I_{w(i+1)})$ by claim 1.
Therefore $\mu_{\mathbb{P}}(I_{wi})\sim\mu_{\mathbb{P}}(I_{w(i+1)})$, and so $p_{k,i}\sim p_{k,i+1}$ by (\ref{eq6}).

\medskip

In the case of $\varepsilon_{k+1}<\delta_{k+1}$, since $n_{k+1}=2$, we have
$$B(x^-(I_{wi2}),\delta_{k+1})\subseteq I_{wi}\subseteq B(x^-(I_{wi2}),3\delta_{k+1})$$
and
$$B(x^+(I_{w(i+1)1}),\delta_{k+1})\subseteq I_{w(i+1)}\subseteq B(x^+(I_{w(i+1)1}),3\delta_{k+1}).$$
Since $\varepsilon_k<\delta_k$, we get
$|x^+(I_{w(i+1)1}-x^-(I_{wi2})|=2\delta_{k+1}+\varepsilon_k\leq5\delta_{k+1}$.
Then $\mu_{\mathbb{P}}(I_{wi})\sim\mu_{\mathbb{P}}(I_{w(i+1)})$ by the doubling property of $\mu_{\mathbb{P}}$, giving
$p_{k,i}\sim p_{k,i+1}$.
This completes the proof of claim 2.

\medskip

{\bf{Claim 3.}} $\mathbb{P}$ is $C$-uniform.

\medskip

Proof of Claim 3.  Since
the $C$-uniformness of $\mathbb{P}$ is equivalent to that every $P_k$ satisfies
\begin{equation}\label{peq1}
p_{k,i+1}+\cdots+p_{k,i+l}\sim p_{k,i+l+1}+\cdots+p_{k,i+2l}
\end{equation}
for all $i\geq0$ and $l\geq1$ with $i+2l\leq n_k$, it suffices to show (\ref{peq1}).
Fix $k$. When $l=1$ or $n_k\leq 3$, the relationship (\ref{peq1}) follows from claim 2 directly.
Next we assume that $l\geq 2$ and $n_k>3$. Let $w\in W_{k-1}$ and $i\geq 0$ be given such that
$i+2l\leq n_k$. Observing that there are balls $B(x,r)$ and $B(y,r)$
with $x,y\in E$ and $|x-y|\leq 3r$ such that
$$E_k\cap B(x,r)\subseteq I_{w(i+1)}\cup\cdots\cup I_{w(i+j)}\subseteq B(x,3r)$$
and
$$E_k\cap B(y,r)\subseteq I_{w(i+j+1)}\cup\cdots\cup I_{w(i+2j)}\subseteq B(y,3r),$$
we immediately get
$$\mu_{\mathbb{P}}(I_{w(i+1)}\cup\cdots\cup I_{w(i+j)})\sim\mu_{\mathbb{P}}(I_{w(i+j+1)}\cup\cdots\cup I_{w(i+2j)}),$$
giving the relationship (\ref{peq1}) by (\ref{eq6}). This completes the proof of claim 3.

\medskip

{\bf{Claim 4.}} $\prod_{j=1}^t{p_{k+j,1}}\sim\prod_{j=1}^t{p_{k+j, n_{k+j}}}$ for all $k\in\Lambda$
and $1\leq t<m_k$.

\medskip

Proof of Claim 4. Let $k\in\Lambda$ be fixed. Without loss of generality, assume that $m_k>1$.
Let $1\leq t<m_k$ be given.
Then $\varepsilon_k<\delta_{k+t}$ by the definition of $m_k$.
Let $w\in W_{k-1}$ and $1\leq i<n_k$ be given. Let $wiu\in W_{k+t}$ be the maximal word
with prefix $wi$ and $w(i+1)v\in W_{k+t}$ be the minimal word with prefix $w(i+1)$.
Then, by (\ref{eq6}) $$\mu_{\mathbb{P}}(I_{wiu})=\mu_{\mathbb{P}}(I_{w})p_{k,i}\prod_{j=1}^t{p_{k+j, n_{k+j}}}$$ and $$\mu_{\mathbb{P}}(I_{w(i+1)v})=\mu_{\mathbb{P}}(I_{w})p_{k,i+1}\prod_{j=1}^t{p_{k+j,1}}.$$
Note that the gap between $I_{wiu}$ and $I_{w(i+1)v}$ is $G_{wi}$. Since $\varepsilon_k<\delta_{k+t}$,
one has $|G_{wi}|<|I_{wiu}|=|I_{w(i+1)v}|$. To complete this proof, we consider two cases:

\medskip

Case 1. $n_{k+t+1}>2$.

\medskip

Arguing as case 2 of claim 2, we have two balls, $B(x,r)$ and $B(y,r)$, with $x,y\in E$ and $|x-y|\leq 6r$, such that
$$B(x,r)\subseteq I_{wiu}\subseteq B(x,4r)\,\,\mbox{and}\,\, B(y,r)\subseteq I_{w(i+1)v}\subseteq B(y,4r),$$
which, combined with the doubling property of $\mu_{\mathbb{P}}$, yields $\mu_{\mathbb{P}}(I_{wiu})\sim\mu_{\mathbb{P}}(I_{w(i+1)v})$.

\medskip

Case 2. $n_{k+t+1}=2$.

\medskip

Consider two subcases: $\varepsilon_{k+t+1}\geq\delta_{k+t+1}$ and $\varepsilon_{k+t+1}<\delta_{k+t+1}$.
In the first subcase, arguing as case 3 of claim 2,
we have two balls, $B(x,r)$ and $B(y,r)$, with $x,y\in E$ and $|x-y|\leq 3r$, such that
$$I_{wiu2}\subseteq E_{k+t}\cap B(x,r)\subseteq I_{wiu}\mbox{ and }I_{w(i+1)v1}\subseteq E_{k+t}\cap B(y,r)\subseteq I_{w(i+1)v}.$$
In the second subcase, arguing as case 3 of claim 2, we have two balls, $B(x,r)$ and $B(y,r)$,
with $x,y\in E$ and $|x-y|\leq 5r$, such that
$$B(x,r)\subseteq I_{wiu}\subseteq B(x,3r)\mbox{ and }B(y,r)\subseteq I_{w(i+1)v}\subseteq B(y,3r).$$
Therefore, in both subcases, we have
$\mu_{\mathbb{P}}(I_{wiu})\sim\mu_{\mathbb{P}}(I_{w(i+1)v})$.

Now $\mu_{\mathbb{P}}(I_{wiu})\sim\mu_{\mathbb{P}}(I_{w(i+1)v})$ is proved.
Since $p_{k,i}\sim p_{k,i+1}$ has been proved in claim 2, we get
$\prod_{j=1}^t{p_{k+j,1}}\sim\prod_{j=1}^t{p_{k+j, n_{k+j}}}$. This completes the proof of claim 4.

\medskip

{\bf{Claim 5.}} $(P_{k+t},P_{k+t})$ is $C$-uniform for all $k\in\Lambda$
and $1\leq t<m_k$.

\medskip

Proof of Claim 5. Since $\mathbb{P}$ has been proved to be $C$-uniform, it suffices to show that
the sum of the first $j$ terms of $P_{k+t}$ is comparable to the sum of the last $j$ terms of $P_{k+t}$
for all $k\in\Lambda$, $1\leq t<m_k$, and $1\leq j\leq n_{k+t}$. Let these integers be given.
Let $w\in W_{k-1}$ and $1\leq i<n_k$ be given. Let $u\in W_{k+t-1}$ be the maximal word
with prefix $wi$ and $v\in W_{k+t-1}$ be the minimal word with prefix $w(i+1)$. Then $G_{wi}$ is
the gap between $I_u$ and $I_v$.
Consider components $I_{un_{k+t}}$, $I_{u(n_{k+t}-1)}$,$\cdots$, $I_{u(n_{k+t}+1-j)}$ and
$I_{v1}$, $I_{v2}$, $\cdots$, $I_{vj}$ of level $k+t$. Since $\varepsilon_k<\delta_{k+t}$,
by the same argument as that of claim 3 we get
$$\mu_{\mathbb{P}}(I_{un_{k+t}}\cup\cdots\cup I_{u(n_{k+t}+1-j)})
\sim\mu_{\mathbb{P}}(I_{v1}\cup\cdots\cup I_{vj}),$$
which, together with (\ref{eq6}), yields
$$p_{k+t,n_{k+t}}+\cdots+p_{k+t,n_{k+t}+1-j}
\sim p_{k+t,1}+\cdots+p_{k+t,j}.$$
This completes the proof of claim 5.

\medskip

{\bf{Claim 6.}} $(P_{k+m_k},P_{k+m_k})$ is $(C,s_k)$-uniform for all $k\in\Lambda$.

\medskip

Proof of Claim 6. Since $\mathbb{P}$ has been proved to be $C$-uniform, it suffices to show that
the sum of the first $j$ terms of $P_{k+m_k}$ is comparable to the sum of the last $j$ terms of $P_{k+m_k}$
for all $k\in\Lambda$ and $s_k\leq j\leq n_{k+m_k}$, where, by the definition, $s_k$ satisfies
$$s_k\delta_{k+m_k}+(s_k-1)\varepsilon_{k+m_k}\leq\varepsilon_k<(s_k+1)\delta_{k+m_k}+s_k\varepsilon_{k+m_k}.$$
Let $k$ and $j$ be given. Let $w\in W_{k-1}$ and $1\leq i<n_k$.
Let $u\in W_{k+m_k-1}$ be the maximal word
with prefix $wi$ and let $v\in W_{k+m_k-1}$ be the minimal word with prefix $w(i+1)$.
Then $G_{wi}$ is the gap between $I_{un_{k+m_k}}$ and $I_{v1}$.

\medskip
Case 1. $s_k=j=1$. In this case, $\delta_{k+m_k}\leq\varepsilon_k\leq2\delta_{k+m_k}+\varepsilon_{k+m_k}$, so we have
$$I_{un_{k+m_k}}\subseteq E\cap B(x^-(G_{wi}),\varepsilon_k)\subseteq I_{un_{k+m_k}}\cup I_{u(n_{k+m_k}-1)}$$
and $$I_{v1}\subseteq E\cap B(x^+(G_{wi}),\varepsilon_k)\subseteq I_{v1}\cup I_{v2}.$$
From claim 2 we get
$$\mu_{\mathbb{P}}(I_{un_{k+m_k}})\sim\mu_{\mathbb{P}}(I_{un_{k+m_k}}\cup I_{u(n_{k+m_k}-1)})\,\,\mbox{and}\,\,\mu_{\mathbb{P}}(I_{v1})\sim\mu_{\mathbb{P}}(I_{v1}\cup I_{v2)}).$$
Since $x^-(G_{wi}), x^+(G_{wi})\in E$ and $|x^-(G_{wi})-x^+(G_{wi})|=\varepsilon_k$, it follows from the doubling
property of $\mu_{\mathbb{P}}$ that
$\mu_{\mathbb{P}}(B(x^-(G_{wi}),\varepsilon_k))\sim\mu_{\mathbb{P}}(B(x^+(G_{wi}),\varepsilon_k))$, so
$$\mu_{\mathbb{P}}(I_{un_{k+m_k}})\sim\mu_{\mathbb{P}}(I_{v1}),$$ which, together with claim 4, yields
$p_{k+m_k, n_{k+m_k}}\sim p_{k+m_k, 1}$.

\medskip
Case 2. $s_k>1$ or $j>1$. In this case, arguing as we did in claim 3, we get
$$\mu_{\mathbb{P}}(I_{un_{k+m_k}}\cup\cdots\cup I_{u(n_{k+m_k}+1-j)})\sim\mu_{\mathbb{P}}(I_{v1}\cup\cdots\cup I_{vj}),$$
which, combined with claim 4, yields $$p_{k+m_k, n_{k+m_k}}+\cdots+p_{k+m_k, n_{k+m_k}+1-j}\sim p_{k+m_k, 1}+\cdots+p_{k+m_k, j}.$$
This completes the proof of claim 6 and the `only if' part is thus proved.

\medskip

\noindent{\bf Proof of the `if' part.} Suppose that $\mathbb{P}$ meets the conditions of Theorem \ref{1}.
We are going to show that $\mu_\mathbb{P}$ is doubling on $E$.

Let $B(x,r)$ be a ball of the real line with $x\in E$ and $r\in(0,1)$.
Let $I$ be the smallest component such that $I\supseteq E\cap B(x,2r)$.
Suppose that $I$ is at the level $k-1$. Then $B(x,2r)$ intersects at least two components of level $k$, so
$2r\geq \varepsilon_k.$
Denote
$$\mathcal{A}=\{L:L\in\mathcal{I}_k,\,\, L\cap B(x,2r)\neq\emptyset\}.$$
Then $\cup_{L\in\mathcal{A}}L\subseteq I$ and $4r\geq(\sharp\mathcal{A}-2)\delta_k+(\sharp\mathcal{A}-1)\varepsilon_k$,
where $\sharp$ denotes the cardinality. We consider four cases as follows.

\medskip

Case 1. $\sharp\mathcal{A}\geq 4$.
\medskip

In this case, we see that $\sharp\mathcal{A}\sim\sharp\{L\in\mathcal{A}: L\subseteq B(x,r)\}$.
Since $P_k$ is $C$-uniform, we get $\mu_{\mathbb{P}}(B(x,r))\sim\mu_{\mathbb{P}}(B(x,2r))$.

\medskip

Case 2. $\sharp \mathcal{A}=3$.

\medskip

In this case, $4r>{\delta_k}$. When $n_{k+1}\geq 4$,
we see that $B(x,r)$ contains at least $[\frac{n_{k+1}}{4}]$
components of level $k+1$ in $L(x)$, where $L(x)\in\mathcal{A}$ is the component
of level $k$ containing $x$. Since $P_{k+1}$ and $P_k$ are $C$-uniform,
we have $$\mu_{\mathbb{P}}(B(x,r))\sim\mu_{\mathbb{P}}(L(x))\sim\mu_{\mathbb{P}}(B(x,2r)).$$
When $n_{k+1}=2$ or $3$, we see from $\sharp \mathcal{A}=3$ that $B(x,r)$ contains at least a
component of level $k+1$. Therefore the last relationship remains true.

\medskip

Case 3. $\sharp\mathcal{A}=2$ and $\varepsilon_{k}\geq\delta_{k}$.

\medskip

In this case, since $2r\geq{\varepsilon_k}$, we have $2r\geq{\delta_k}$. By the same argument
as that of Case 2 we get $\mu_{\mathbb{P}}(B(x,r))\sim\mu_{\mathbb{P}}(B(x,2r))$.

\medskip

Case 4. $\sharp\mathcal{A}=2$ and $\varepsilon_{k}<\delta_{k}$.

\medskip

Let $I_{wi}$ and $I_{w(i+1)}$ be members of $\mathcal{A}$, where $w\in W_{k-1}$ and $1\leq i<n_k$.
Without loss of generality, assume that $x\in I_{wi}$. Let $J(x)$ be the biggest component such that
$x\in J(x)\subseteq B(x,r)$. Then $J(x)=I$ or $J(x)\subseteq I_{wi}$. It is obvious that
$\mu_{\mathbb{P}}(B(x,r))=\mu_{\mathbb{P}}(B(x,2r))$ for $J(x)=I$. Also, if $J(x)=I_{wi}$,
we easily see from the $C$-uniformness of $P_k$ that
$$\mu_{\mathbb{P}}(B(x,r))\geq\mu_{\mathbb{P}}(I_{wi})\geq C\mu_{\mathbb{P}}(I_{wi}\cup I_{w(i+1)})
\geq C\mu_{\mathbb{P}}(B(x,2r)).$$
Next assume $J(x)\neq I_{wi}$. Thus $J(x)$ is at the level $k+t$ for some $t\geq1$.
Let $I_u$ be the component of level $k+t-1$ containing $J(x)$, where $u\in W_{k+t-1}$.
Then we have $B(x,r)\subset I_u\cup G\cup I_v$ and $r<\delta_{k+t-1}$ by the maximality of $J(x)$,
where $I_v$ is the nearest component of level $k+t-1$ on the right of $I_u$
and $G$ is the gap between $I_u$ and $I_v$. Therefore, $B(x,2r)$ meets at most four components of level $k+t-1$,
that is, $\sharp\mathcal{B}\leq 4$, where
$$\mathcal{B}=\{L:L\mbox{ is a component of $E_{k+t-1}$ and}\,\, L\cap B(x,2r)\neq\emptyset\}.$$
Clearly, $I_u\in\mathcal{B}$. We claim that $\mu_{\mathbb{P}}(L)\sim\mu_{\mathbb{P}}(I_u)$ for all $L\in\mathcal{B}$.
In fact, let $L$ and $K$ be two adjacent components in $\mathcal{B}$,
then the gap between $L$ and $K$ is at the level $k+j$ for some $0\leq j\leq t-1$.
When the gap is at the level $k+t-1$, we immediately get
$\mu_{\mathbb{P}}(L)\sim\mu_{\mathbb{P}}(K)$ from the $C$-uniformness of $P_{k+t-1}$.
The gap being at the level $k+j$ for some $1\leq j\leq t-2$ is possible only when
$$\varepsilon_{k+j}\leq2r<2\delta_{k+t-1}<\delta_{k+j},$$  and if it happened, we would have $m_{k+j}\geq t-1$. 
If $m_{k+j}= t-1$, then $s_{k+j}=1$, so the $(C, s_{k+j})$-uniformness
of $(P_{k+j+m_{k+j}},P_{k+j+m_{k+j}})$ implies $\mu_{\mathbb{P}}(L)\sim\mu_{\mathbb{P}}(K)$.
If $m_{k+j}> t-1$, we may use the condition (\ref{eq11}) to
get $\mu_{\mathbb{P}}(L)\sim\mu_{\mathbb{P}}(K)$. Finally, we consider the case where the gap is at the level $k$.
Note that $t\leq m_k+1$, in fact, if not, we would get from the definition of $m_k$
that $2r\leq2\delta_{k+m_k+1}<\delta_{k+m_k}\leq\varepsilon_k$,
contradicting $2r\geq\varepsilon_k$. Also, if $t=m_k+1$ happened, $(P_{k+m_k},P_{k+m_k})$ would have been $C$-uniform. In fact,
if $t=m_k+1$ then $\delta_{k+m_k}\leq{\varepsilon_k}\leq 2r\leq2\delta_{k+m_k}$
by the definition of $m_k$, so $s_k=1$, and hence $(P_{k+m_k},P_{k+m_k})$ is $C$-uniform,
because it has been assumed to be $(C, s_k)$-uniform. Then $\mu_{\mathbb{P}}(L)\sim\mu_{\mathbb{P}}(K)$
follows from the above discussion and the assumptions of theorem \ref{1}. This proves the claim.

Now let
$$\mathcal{C}=\{L:L\mbox{ is a component of $E_{k+t}$ and}\,\, L\subseteq I_u\cap B(x,r)\}.$$
Two possible subcases may happen.

\medskip
Subcase 1. $B(x,2r)\supseteq I_u$.
\medskip

In this subcase, we have
$\sharp\mathcal{C}\geq[\frac{n_{k+t}}{2}]$. It follows from the $C$-uniformness of $P_{k+t}$ that
$$\mu_{\mathbb{P}}(B(x,r))\geq C\mu_{\mathbb{P}}(I_u).$$
Since $\sharp\mathcal{B}\leq 4$, we get from the above claim that
$\mu_{\mathbb{P}}(B(x,r))\sim\mu_{\mathbb{P}}(B(x,2r))$.

\medskip
Subcase 2. $B(x,2r)\subset I_u\cup G\cup I_v$.
\medskip

In this subcase, it is clear that $2r<|I_u|=\delta_{k+t-1}$,
so we have $t\leq m_k$ by the definition of $m_k$. Observing that $I_u\subseteq I_{wi}$
and $I_v\subseteq I_{w(i+1)}$, with $u\in W_{k+t-1}$ being the maximal word with
prefix $wi$ and $v\in W_{k+t-1}$ being the minimal word with
prefix $w(i+1)$, we see that $G$ is a gap of level $k$.
When $1\leq t<m_k$, we use the assumption (a) to get $\mu_{\mathbb{P}}(B(x,r))\sim\mu_{\mathbb{P}}(B(x,2r))$.
When $t=m_k$, remembering that $s_k$ has been defined to satisfy
$$s_k\delta_{k+m_k}+(s_k-1)\varepsilon_{k+m_k}\leq\varepsilon_k<(s_k+1)\delta_{k+m_k}+s_k\varepsilon_{k+m_k},$$
which, together with $2r\geq\varepsilon_{k}$, yields
$2r\geq s_k\delta_{k+m_k}+(s_k-1)\varepsilon_{k+m_k}$, we may use the $C$-uniformness of $P_{k+m_k}$,
the $(C, s_k)$-uniformness of $(P_{k+m_k},P_{k+m_k})$, and the assumption (\ref{eq11}) to get
$\mu_{\mathbb{P}}(B(x,r))\sim\mu_{\mathbb{P}}(B(x,2r))$.
This completes the proof of the `if' part.

\section{Proof of Theorem \ref{t11}}


The proof of theorem 2 may go in step as that of theorem 1. For the proof of the `only if' part, suppose
that $\mu$ is doubling on the uniform Cantor set $E$.
We shall prove that the sequence $\{P_w\}_{w\in W}$ satisfies the conditions of theorem 2.
As we did in theorem 1, the proof consists of the following claims.

\medskip

{\bf{Claim A.}} $p_{w2}\leq Cp_{w1}$ and $p_{w(n_k-1)}\leq Cp_{wn_k}$ for all $w\in W_{k-1}$ and $k\geq 1$.

\medskip

{\bf{Claim B.}} $p_{wi}\sim p_{w(i+1)}$
for all $w\in W_{k-1}$, $1\leq i<n_k$, and $k\geq 1$.

\medskip

{\bf{Claim C.}} $P_w$ is $C$-uniform for all ${w\in W}$.

\medskip

{\bf{Claim D.}} $\prod_{j=1}^{t}{p_{w(i+1)1^j}}\sim\prod_{j=1}^{t}p_{wiu_j}$ for all $w\in W_{k-1}$,
$1\leq i<n_k$, $1\leq t<m_k$, and $k\in \Lambda$.

\medskip

{\bf{Claim E.}} $(P_{wiu_{t-1}},P_{w(i+1)1^{t-1}})$ is $C$-uniform for all $w\in W_{k-1}$, $1\leq i<n_k$,
$1\leq t<m_k$, and $k\in \Lambda$.

\medskip

{\bf{Claim F.}} $(P_{wiu_{m_k-1}},P_{w(i+1)1^{m_k-1}})$ is $(C,s_k)$-uniform for all $k\in\Lambda$.

\medskip

If, for every fixed integer $k\geq 1$, the vectors $P_w$ are the same for all $w\in W_{k-1}$,
then the above claims are exactly those in theorem 1. Under this condition, {\it claim A} has a
simpler proof, because the word $w\in W_{k-1}$ in question may be assumed to be the minimal or
maximal words. Now, without assuming this condition, we are going to prove {\it claim A}. Instead,
we will show the following

\begin{lemma}\label{m}
$p_{w2}\sim p_{w1}$ and $p_{w(n_k-1)}\sim p_{wn_k}$ for all $w\in W_{k-1}$ and $k\geq 1$.
\end{lemma}

Proof of Lemma \ref{m}. We only prove $p_{w2}\sim p_{w1}$. The proof of the other relationship
is obviously similar.
Let $k\geq 1$ and $w\in W_{k-1}$ be given. To prove $p_{w2}\sim p_{w1}$, it suffices to show
$\mu(I_{w2})\sim\mu(I_{w1})$. Without loss of generality, assume that
$w\in W_{k-1}$ is neither minimal nor maximal.

\medskip

Case 1. $\varepsilon_k\geq\delta_k$.

\medskip

Let $I_u$ be the component of level $k$ next to the left of $I_{w1}$. Let $I_v$ be the component of
level $k$ next to the right of $I_{w2}$.
Let $G$ be the gap
between $I_u$ and $I_{w1}$. Let $G'$ be the gap between $I_{w2}$ and $I_v$. Then $G'$ is a gap
of level $k$ if and only if $n_k>2$. Remember
that $\mu$ is doubling on $E$ and that $\varepsilon_k\geq\delta_k$ has been assumed.

Noting that, if $\delta_k\leq|G|<\varepsilon_k$ then
$$I_{w1}=E_k\cap B(x^+(I_{w1}),|G|)\mbox{ and } I_{u}=E_k\cap B(x^-(I_{u}),|G|),$$
and that, if $|G|<\delta_k$ then
$$I_{w1}=E_k\cap B(x^+(I_{w1}),\delta_k)\mbox{ and }I_{u}=E_k\cap B(x^-(I_{u}),\delta_k),$$
we find that, if $|G|<\varepsilon_k$, then
\begin{equation}\label{rr}
\mu(I_{w1})\sim \mu(I_u).
\end{equation}
Arguing as we just did, if $|G'|<\varepsilon_k$, then
\begin{equation}\label{rrr}
\mu(I_{w2})\sim \mu(I_v).
\end{equation}

Now we apply (\ref{rr}) and (\ref{rrr}) to prove $\mu(I_{w1})\sim\mu(I_{w2})$. We consider nine subcases.

\medskip

Subcase 1. $|G|\geq\varepsilon_k$ and $|G'|\geq\varepsilon_k$.

Since $\varepsilon_k\geq\delta_k$ has been assumed,
we have $$I_{w1}=E_k\cap B(x^+(I_{w1}),\varepsilon_k)\,\,\mbox{and}\,\, I_{w2}=E_k\cap B(x^-(I_{w2}),\varepsilon_k),$$ and so
$\mu(I_{w1})\sim\mu(I_{w2}).$

\medskip

Subcase 2. $|G|<\delta_k$ and $|G'|<\delta_k$.

In view of (\ref{rr}) and (\ref{rrr}),
$\mu(I_{u})\sim\mu(I_{w1})$
and $\mu(I_{v})\sim\mu(I_{w2})$.
Since $|G'|<\delta_k\leq\varepsilon_k$, we have $n_k=2$, as we have said.
We see that the gap next to the right of $I_v$ belongs to $\mathcal{G}_k$. Then we have
$$
\mu(I_{w1})+\mu(I_{u})=\mu(B(x^+(I_{u}),\varepsilon_k+\delta_k))
$$
and
$$
\mu(I_{w2})+\mu(I_{v})=\mu(B(x^-(I_{v}),\varepsilon_k+\delta_k)),
$$
and so $\mu(I_{w1})\sim\mu(I_{w2})$.

\medskip

Subcase 3. $\delta_k\leq|G|<\varepsilon_k$ and $|G'|\geq\varepsilon_k$.

In view of (\ref{rr}), $\mu(I_{w1})\sim\mu(I_{u})$. Observing that
$$I_{w1}\subseteq E_k\cap B(x^+(I_{w1}),\varepsilon_k)\subseteq I_u\cup I_{w1}\,\,\mbox{and}\,
\,I_{w2}=E_k\cap B(x^-(I_{w2}),\varepsilon_k),$$
we get $\mu(I_{w1})\sim\mu(I_{w2})$.

\medskip

Subcase 4. $\delta_k\leq|G'|<\varepsilon_k$ and $|G|\geq\varepsilon_k$.

Arguing as we just did in subcase 3, we have $\mu(I_{w1})\sim\mu(I_{w2})$.

\medskip

Subcase 5. $\delta_k\leq|G|<\varepsilon_k$ and $|G'|<\delta_k$.

We have
$\mu(I_{u})\sim\mu(I_{w1})$ and $\mu(I_{v})\sim\mu(I_{w2})$. We also have
$$I_{w1}\subseteq E_k\cap B(x^-(I_{w1}),\varepsilon_k+\delta_k)\subseteq I_u\cup I_{w1}$$
and
$$\mu(I_{w2})+\mu(I_{v})=\mu(B(x^-(I_{v}),\varepsilon_k+\delta_k)),$$
and so
$\mu(I_{w1})\sim\mu(I_{w2}).$

\medskip

Subcase 6. $\delta_k\leq|G'|<\varepsilon_k$ and $|G|<\delta_k$.

We get $\mu(I_{w1})\sim\mu(I_{w2})$ by an argument analogous to case 5.

\medskip

Subcase 7. $\delta_k\leq|G|<\varepsilon_k$ and $\delta_k\leq|G'|<\varepsilon_k$.

We have
$\mu(I_{u})\sim\mu(I_{w1})$ and $\mu(I_{v})\sim\mu(I_{w2})$. We also have
$$I_{w1}\subseteq E_k\cap B(x^-(I_{w1}),\varepsilon_k+\delta_k)\subseteq I_u\cup I_{w1}$$
and
$$I_{w2}\subseteq E_k\cap B(x^+(I_{w2}),\varepsilon_k+\delta_k)\subseteq I_v\cup I_{w2},$$
giving $\mu(I_{w1})\sim\mu(I_{w2}).$

\medskip

Subcase 8. $|G|<\delta_k$ and $|G'|\geq\varepsilon_k$.

We have
$\mu(I_{u})\sim\mu(I_{w1})$. Noting that
$$I_{w2}=E_k\cap B(x^-(I_{w2}),\varepsilon_k)\,\,\mbox{and}\,\,I_{w1}\subseteq E_k\cap B(x^+(I_{w1}),
\varepsilon_k)\subseteq I_u\cup I_{w1},$$
we get $\mu(I_{w1})\sim\mu(I_{w2})$.

\medskip

Subcase 9. $|G'|<\delta_k$ and $|G|\geq\varepsilon_k$.

Arguing as in case 8, we have $\mu(I_{w1})\sim\mu(I_{w2})$.

\medskip
Case 2. $\varepsilon_k<\delta_k$, $n_{k+1}=2$, and $\varepsilon_{k+1}\geq\delta_{k+1}$.
\medskip

We consider two subcases: $\varepsilon_k<\delta_{k+1}$ and $\varepsilon_k\geq\delta_{k+1}$. When $\varepsilon_k<\delta_{k+1}$,
one has $\mu(I_{w12})\sim\mu(I_{w21})$ by comparing the sizes of
$B(x^-(I_{w12}),\delta_{k+1})$ and $B(x^+(I_{w21}),\delta_{k+1})$ in measure $\mu$, which, combined with an observation that  $I_{w1}\subset B(x^-(I_{w12}),2\varepsilon_{k+1})$, yields $$\mu(I_{w1})\leq C\mu(B(x^-(I_{w12}),\varepsilon_{k+1}))\leq C(\mu(I_{w12})+\mu(I_{w21}))\leq C\mu(I_{w2}).$$
As $\mu(I_{w2})\leq C\mu(I_{w1})$ may be obtained similarly, we have $\mu(I_{w1})\sim\mu(I_{w2})$. When $\varepsilon_k\geq\delta_{k+1}$, since $3\varepsilon_{k+1}\geq\delta_k>\varepsilon_k$ by the assumptions of case 2,
one also has $\mu(I_{w12})\sim\mu(I_{w21})$ by comparing
$B(x^-(I_{w12}),r_{k})$ and $B(x^+(I_{w21}),r_{k})$, where $r_{k}=\min\{\varepsilon_k,\varepsilon_{k+1}\}$. Therefore, arguing as in the previous subcase gives
$\mu(I_{w1})\sim\mu(I_{w2})$.

\medskip
Case 3. $\varepsilon_k<\delta_k$, $n_{k+1}=2$, and $\varepsilon_{k+1}<\delta_{k+1}$.
\medskip

Case 4. $\varepsilon_k<\delta_k$ and $n_{k+1}>2$.

\medskip

For cases 3 and 4, the proof of $\mu(I_{w1})\sim\mu(I_{w2})$
is the same as that of claim 2 in section 2. We omit it. $\hfill\Box$

\medskip

Based on lemma \ref{m}, {\it claims B - F} can be proved in the same way as {\it claims 2 - 6}
of section 2. This completes the proof of the `only if' part.

\medskip

For the proof of the `if' part, suppose that the sequence $\{P_w\}_{w\in W}$ satisfies the conditions of theorem 2.
We have to prove that $\mu$ is doubling on the uniform Cantor set $E$. Arguing as we did in the proof of theorem 1,
the doubling property of $\mu$ can be proved similarly, and so it is omitted here.

\section{Proof of Theorem \ref{t2}}

Let $E=E(\mathbf{n},\mathbf{c},\Lambda, \{m_k\}_{k\in\Lambda}, \{s_k\}_{k\in\Lambda})$ be a uniform Cantor set.
The proof of Theorem \ref{t2} is based on the following lemmas.

\begin{lemma}\label{lem1}
If $E$ carries a doubling measure which can be extended to a doubling measure on $[0,1]$, then
$\{n_kc_k\}_{k=1}^\infty\in\cup_{q\geq1}\ell^q$.
\end{lemma}
\noindent{\bf Proof.\,} If $\{n_kc_k\}_{k=1}^\infty\not\in\cup_{q\geq1}\ell^q$, then we have $\nu(E)=0$
for all doubling measures $\nu$ on $[0,1]$ (see \cite{hww}), a contradiction.
$\hfill\Box$

\begin{lemma}\label{lem2}
Suppose that $E$ carries a doubling measure which can be extended to a doubling measure
on $[0,1]$ and that $\sup_{k\geq1}n_k<\infty$. Then $\sup_{k\in\Lambda}m_k=\infty$.
\end{lemma}
\noindent{\bf Proof.\,} If $\sup_{k\in\Lambda}m_k<\infty$, then there is an integer $N$ such
that $\varepsilon_k\geq\delta_{k+N}$ for all $k\geq 1$, i.e.
$$\prod_{i=k}^{k+N}\frac{1-(n_i-1)c_i}{n_i}\leq c_k.$$
It follows that
\begin{equation}\label{jj}
(\frac{1-\sup _{i\geq1}(n_i-1)c_i}{\sup _{i\geq1}n_i})^{N+1}\leq c_k
\end{equation}
for all $k\geq 1$. Since $E$ carries a doubling measure which can be extended
to a doubling measure
on $[0,1]$, we have from Lemma \ref{lem1} that $\{n_kc_k\}_{k=1}^\infty\in\cup_{q\geq1}\ell^q$.
This implies $\sup_{k\geq 1}(n_k-1)c_k<1$, which, together with (\ref{jj}) and the assumption
$\sup_{k\geq1}n_k<\infty$, gives $\inf_{k\geq 1}c_k>0$. This implies that $E$ is porous, so
we have $\nu(E)=0$ for all doubling measures $\nu$ on $[0,1]$ (see \cite{Wu1}), a contradiction. $\hfill\Box$

\medskip

\noindent{\bf Proof of Theorem \ref{t2}.} Suppose $\{n_kc_k\}_{k=1}^\infty\in\ell^1$.
Then $\mathcal{L}(E)>0$ by (\ref{eq4}). Since $\mathbb{P}$ is an ultimately $1$-uniform $\mathbf{n}$-matching
probability sequence, there is an integer $k_0\geq 1$ such that
$p_{k,i}=\frac{1}{n_k}$ for all $k\geq k_0$ and $1\leq i\leq n_k$. By this, we have
$$\mu_{\mathbb{P}}=\sum_{I\in\mathcal{I}_{k_0-1}}\frac{\mu_{\mathbb{P}}(I)}{\mathcal{L}(E\cap I)}\mathcal{L}|_{E\cap I}.$$
Let
$$\nu=\sum_{I\in\mathcal{I}_{k_0-1}}\frac{\mu_{\mathbb{P}}(I)}{\mathcal{L}(E\cap I)}\mathcal{L}|_{I}
+\sum_{J\in\mathcal{G}_{k_0-1}}\mathcal{L}|_J.$$
Since $\mathcal{I}_{k_0-1}\cup\mathcal{G}_{k_0-1}$ forms a partition of $[0,1]$ by finite intervals
and the restrictions of the measure $\nu$ to these intervals are comparable to Lebesgue measure, we easily see that $\nu$ is doubling on $[0,1]$.
Moreover, since $E\cap(\cup_{G\in\mathcal{G}_{k_0-1}}G)=\emptyset$, we have
$\mu_{\mathbb{P}}=\nu|_E.$  This proves that $\mu_{\mathbb{P}}$ can be extended to a doubling measure on $[0,1]$.

Conversely, suppose that $\mu_{\mathbb{P}}$ can be extended to a doubling measure $\nu$
on $[0,1]$.
Then we have from Lemma \ref{lem1} that
$\{n_kc_k\}_{k=1}^\infty\in\ell^q$ for some $q\geq1$,
which yields $\lim_{k\rightarrow\infty}n_kc_k=0$.
We are going to show $\{n_kc_k\}_{k=1}^\infty\in\ell^1$.

Since $\lim_{k\rightarrow\infty}n_kc_k=0$ and $\mathbb{P}$ is ultimately $1$-uniform,
we may choose an integer $k_0\geq 1$ such that
\begin{equation}\label{jia2}
\mbox{$n_kc_k<1/3$\, for all $k\geq k_0$}
\end{equation}
and that
\begin{equation}\label{jia3}
\mbox{$p_{k,i}=\frac{1}{n_k}$ for all $k\geq k_0$ and $1\leq i\leq n_k$.}
\end{equation}
Therefore $\{k:k\geq k_0\}\subseteq\Lambda$, where $\Lambda$ is as in (\ref{eq7}).
Let $k\geq k_0$ be fixed. By the definitions of integers $m_k$ and $s_k$, we have
\begin{equation}\label{jia1}
s_k\delta_{k+m_k}+(s_k-1)\varepsilon_{k+m_k}\leq\varepsilon_k<(s_k+1)\delta_{k+m_k}+s_k\varepsilon_{k+m_k}.
\end{equation}
Since $\varepsilon_k=c_k\delta_{k-1}$ by the construction of the uniform Cantor set $E$, we see from (\ref{eq3}) that
the right-hand inequality of (\ref{jia1}) can be rewritten as
$$c_k<\frac{(s_k+1)(1-(n_{k+m_k}-1)c_{k+m_k})+s_kn_{k+m_k}c_{k+m_k}}{n_{k+m_k}}\prod_{i=k}^{k+m_k-1}\frac{1-(n_i-1)c_i}{n_i}.$$
Therefore
\begin{equation}\label{p21}
n_kc_k\leq\frac{3s_k}{n_{k+1}n_{k+2}\cdots n_{k+m_k}}.
\end{equation}
On the other hand, let $w\in W_{k-1}$, $1\leq i<n_k$, and let $G_{wi}\in\mathcal{G}_k$ be
the corresponding gap of level $k$. Then $|G_{wi}|=\varepsilon_k$. Let $J=[x^-(G_{wi})-\varepsilon_k,x^-(G_{wi})]$ be an interval,
where $x^-(G_{wi})$ is the left endpoint. We see from the left-hand inequality of (\ref{jia1}) that
$J$ contains at least $s_k$ components of level $k+m_k$ in $I_{wi}$. Since $\nu$ is an extension
of $\mu_{\mathbb{P}}$ and $\nu$ is doubling on $[0,1]$, we get from (\ref{jia3}) that
$$\nu(G_{wi})\geq C\nu(J)\geq\frac{Cs_k\mu_{\mathbb{P}}(I_{wi})}{n_{k+1}\cdots n_{k+m_k}}.$$
Now, summing over all $1\leq i<n_k$ yields
\begin{equation}\label{gp}
\nu(\bigcup_{i=1}^{n_{k}-1}G_{wi})\geq\frac{Cs_k}{n_{k+1}\cdots n_{k+m_k}}\mu_{\mathbb{P}}(\bigcup_{i=1}^{n_{k}-1}I_{wi}).
\end{equation}
Since ${\mathbb{P}}$ is ultimately $1$-uniform, we have $\mu_{\mathbb{P}}(\cup_{i=1}^{n_{k}-1}I_{wi})\sim\mu_{\mathbb{P}}(I_w)$.
Summing over all $w\in W_{k-1}$, we get from (\ref{p21}) that
$$
\nu(\bigcup_{J\in\mathcal{G}_k}J)\geq\frac{Cs_k}{n_{k+1}\cdots n_{k+m_k}}
\geq Cn_kc_k
$$
Finally, summing over all $k\geq k_0$ yields
$$\infty>\nu([0,1])\geq \sum_{k=k_0}^\infty\nu(\bigcup_{J\in\mathcal{G}_k}J)
\geq C\sum_{k=k_0}^\infty n_kc_k$$
This proves $\{n_kc_k\}\in\ell^1$.

\medskip

\noindent{\bf Acknowledgement}\, The authors wish to thank the referee for
their valuable comments.

\noindent Chun WEI and Zhixiong WEN  \\
Department of Mathematics \\
Huazhong University of Science and Technology \\ Wuhan 430074, China
\\ \\ Shengyou WEN\\ Department of Mathematics\\
Hubei University \\ Wuhan 430062, China \\
E-mail address: sywen\_65@163.com

\end{document}